\documentclass{amsart}
\pdfoutput=1

\usepackage[utf8]{inputenc}

\usepackage{amsmath}
\usepackage{amsfonts}
\usepackage{amssymb}
\usepackage{amsthm}
\usepackage{enumerate}
\usepackage{enumitem}
\usepackage{color}
\usepackage{stmaryrd}
\usepackage{slashed}
\usepackage{caption}

\usepackage{setspace}
\usepackage{mathtools}

\usepackage{relsize}
\usepackage{tikz-cd}
\usepackage{tikz}
\usetikzlibrary{decorations.markings}
\usetikzlibrary{decorations.pathmorphing}
\usetikzlibrary{cd}
\usepackage{comment}
\usepackage{wrapfig}

\usepackage{geometry}
\usepackage[colorlinks=true,linkcolor=green!50!black,citecolor=red!75!black]{hyperref}
\usepackage{tikz-3dplot}
\usepackage{bbm}
\usepackage[T1]{fontenc}

\newcommand{\Bord}{\operatorname{Bord}}

\newcommand{\Aut}{\operatorname{Aut}}
\newcommand{\Hom}{\operatorname{Hom}}

\newcommand{\End}{\operatorname{End}}
\newcommand{\pt}{\operatorname{pt}}
\newcommand{\Mod}{\operatorname{Mod}}

\newcommand{\Map}{\operatorname{Map}}

\newcommand{\Ab}{\operatorname{Ab}}

\newcommand{\Fun}{\operatorname{Fun}}

\newcommand{\Pic}{\operatorname{Pic}}

\newcommand{\Z}{\mathbb{Z}}

\newcommand{\R}{\mathbb{R}}
\newcommand{\E}{\mathbb{E}}

\newcommand{\Cat}{\operatorname{Cat}}

\newcommand{\Sp}{\operatorname{Sp}}
\newcommand{\colim}{\operatorname{colim}}
\newcommand{\GL}{\mathrm{GL}}
\newcommand{\gl}{\mathrm{gl}}
\newcommand{\fr}{\mathrm{fr}}
\newcommand{\fd}{\mathrm{fd}}
\newcommand{\TFT}{\mathrm{TFT}}

\newtheorem{theorem}{Theorem}[section]

\newtheorem{proposition}[theorem]{Proposition}
\newtheorem{lemma}[theorem]{Lemma}

\newtheorem{corollary}[theorem]{Corollary}

\theoremstyle{definition}

\newtheorem{definition}[theorem]{Definition}

\theoremstyle{remark}

\newtheorem{remark}[theorem]{Remark}
\newtheorem{example}[theorem]{Example}

\definecolor{Blue} {rgb} {0.282352,0.239215,0.803921}
\definecolor{Green} {rgb} {0.133333,0.545098,0.133333}

\newcounter{jfc}

\title{The $J$-homomorphism from localising the cobordism hypothesis}

\author[]{Amartya Shekhar Dubey}
\address{Amartya Shekhar Dubey, Department of Mathematical Sciences, Norwegian University of Science and Technology, Trondheim, Norway}
\email{amartya.dubey@ntnu.no}

\author[]{Luuk Stehouwer}
\address{Luuk Stehouwer, Durham University, Department of Mathematical Sciences, Upper Mountjoy, Durham DH1 3LE, UK}
\email{stehouwer.d.luuk@durham.ac.uk}

\begin{document}

\begin{abstract}
Lurie's cobordism hypothesis classifies fully extended framed topological field theories valued in a symmetric monoidal $(\infty,n)$-category $\mathcal{C}$ with duals by the $\infty$-groupoid $\mathcal{C}^{\mathrm{fd}}$ of fully dualisable objects. This classification results in an $O(n)$-action on $\mathcal{C}^{\mathrm{fd}}$, where $O(n)$ acts on the framed bordism category by change of framing. Lurie observed that for a connective spectrum, this action is ``a version of the classical $J$-homomorphism.'' We make this identification explicit: with respect to the right action of $O(n)$ on framings, the $O(n)$-action on the localisation $L(\mathrm{Bord}_n^{\mathrm{fr}}) \simeq S$ is precisely $\ominus J_n$, the inverse of the unstable $J$-homomorphism $J_n \colon O(n) \to \mathrm{GL}_1(S)$. The identification of the action with the unstable $J$-homomorphism in the setting of smooth stacks is due to Grady--Pavlov \cite[Theorem 5.0.6]{grady2021geometric}; we give an independent proof in Lurie's original formulation of the cobordism hypothesis, using the abstract Pontryagin--Thom construction. The sign (i.e.\ whether the action is $J_n$ or $\ominus J_n$) depends on the convention for the $O(n)$-action on framings, and is relevant for $\mathrm{pin}^\pm$ theories and hence for quantum systems with time-reversal symmetry. We compute the induced action on homotopy groups of the sphere spectrum in low degrees, yielding explicit invariants for invertible TFTs in dimensions $1$, $2$, and $3$, and connect these to the Kitaev chain, the framing anomaly, the spin-statistics theorem, and anomaly inflow. We recover the Galatius--Madsen--Tillmann--Weiss theorem as a consequence of the cobordism hypothesis and relate the homotopy fixed-point and mapping-spectrum perspectives on invertible TFTs via the duality between tangential and normal framings.

\end{abstract}

\maketitle
\begin{figure*}[!htbp]
\centering
\includegraphics[width=0.68\textwidth]{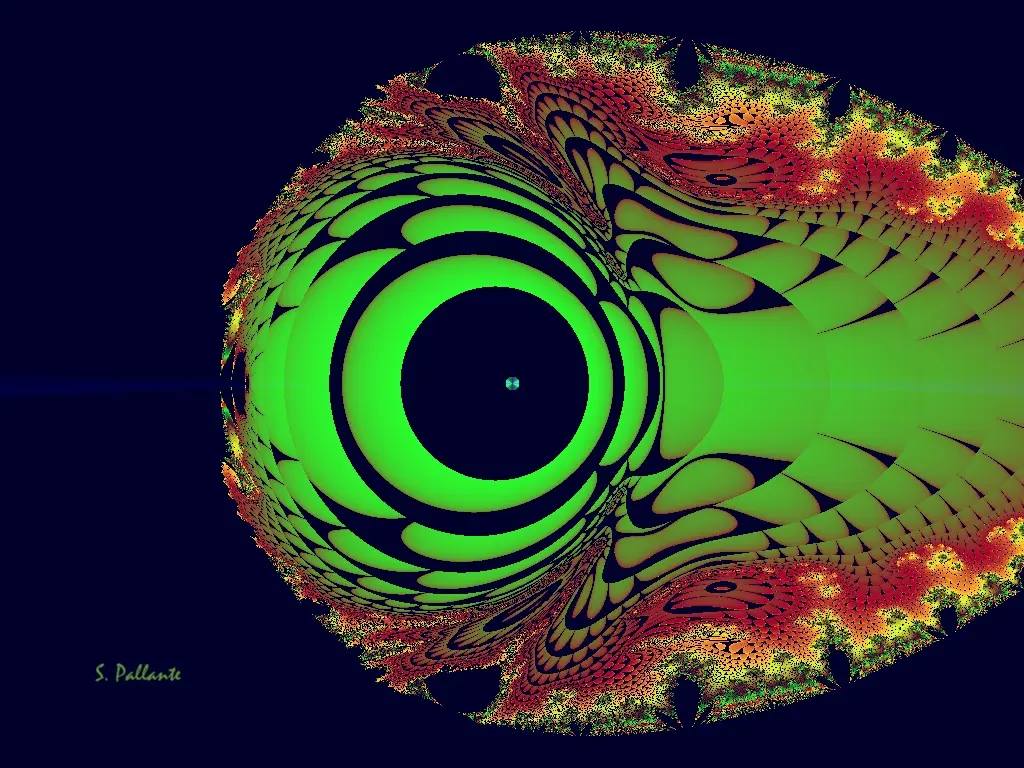}
\caption*{\emph{Green Anguish}! Art by Sergio Pallante \cite{pallante2026artwork}}
\end{figure*}
\clearpage

\section{Introduction}
The cobordism hypothesis, stated by Baez--Dolan \cite{baez1995higher} and sketched by Lurie \cite{MR2555928}, is a fundamental result in the theory of topological quantum field theories (TFTs). It asserts that a fully extended framed $n$-dimensional TFT valued in a symmetric monoidal $(\infty,n)$-category $\mathcal{C}$ with duals is completely determined by its value on the positively framed point. More precisely, the space $\TFT_n^{\fr}(\mathcal{C})$ of such theories is equivalent to the space\footnote{Throughout this paper, we use the algebraic topologist's convention that the word `space' refers to `homotopy type up to higher coherent homotopy', also known as $\infty$-groupoid or anima.} $\mathcal{C}^{\fd}$ of fully dualisable objects of $\mathcal{C}$. Since $O(n)$ acts on $\Bord_n^{\fr}$ by change of framing this induces an $O(n)$-action on $\mathcal{C}^{\fd}$.

A central consequence of the cobordism hypothesis is the classification of \emph{invertible} TFTs, which are those valued in a connective spectrum $X$, such as the sphere spectrum $S$ or the maximal Picard subgroupoid of some symmetric monoidal $(\infty,n)$-category. Invertible TFTs describe the low-energy limits of symmetry-protected topological (SPT) phases and the anomalies of quantum field theories~\cite{freed2021reflection,Witten_2016}. The classification of such theories reduces, via the cobordism hypothesis, to understanding the $O(n)$-action on the localisation $L(\Bord_n^{\fr})$ of the framed bordism category, which is a connective spectrum. The Galatius--Madsen--Tillmann--Weiss (GMTW) theorem \cite{MR2506750} identifies $L(\Bord_n^{B})$ with the Madsen--Tillmann spectrum $MTB$ for any tangential structure $B \to BO(n)$, providing a powerful computational tool.

Lurie observed that the $O(n)$-action on the space of fully dualisable objects of a Picard $\infty$-groupoid is ``a version of the classical $J$-homomorphism'' from stable homotopy theory. The $J$-homomorphism $J_n\colon O(n) \to \GL_1(S)$, defined by the one-point compactification of linear isometries, is one of the most studied maps in algebraic topology; its image in the stable homotopy groups of spheres is the subject of the Adams conjecture \cite{adams1963groups} and controls the chromatic height-$1$ layer of the sphere spectrum. This identification has also been established by Grady-Pavlov \cite[Theorem~5.0.6]{grady2021geometric} in the framework of geometric cobordism categories.

In this paper, we make Lurie's observation completely explicit by giving a proof in Lurie's original formulation of the cobordism hypothesis. Our main theorem identifies the $O(n)$-action on $L(\Bord_n^{\fr}) \simeq S$ with the inverse of the unstable $J$-homomorphism, $\ominus J_n$ (Theorem~\ref{thm:main}), using the abstract Pontryagin--Thom collapse map. The sign depends on the convention for the $O(n)$-action on framings; with the right-action convention used in this paper, the action is $\ominus J_n$. An equivalent result has been obtained by Grady-Pavlov \cite[Theorem 5.0.6]{grady2021geometric} in the framework of geometric cobordism categories; our proof works directly within Lurie's formulation. We explain how the choice of convention (right vs.\ left action of $O(n)$ on framings) determines whether the tangential action is $\ominus J_n$ or $J_n$, and how this relates to the normal action appearing in the construction of Madsen--Tillmann spectra. This distinction, while invisible in low-degree homotopy groups, is relevant for the classification of pin$^\pm$ theories and the physics of time-reversal symmetry.

With the action identified, we compute its effect on the homotopy groups of $S$ in low degrees (Proposition~\ref{prop:action-summary}), yielding explicit invariants for invertible TFTs in dimensions $n = 1, 2, 3$. We then connect these computations to the physics of the Kitaev chain, the framing anomaly of chiral conformal field theories, the spin-statistics theorem, and anomaly inflow \cite{witten1989quantum, witten1992holomorphic, freed2021reflection}. Finally, we recover the GMTW theorem as a formal consequence of the cobordism hypothesis, and discuss the relationship between the homotopy fixed point spectral sequence and the Adams spectral sequence for $MTB$.

Our approach here offers two perspectives on the classification of invertible TFTs: the homotopy fixed point formulation $X^{hB}$ (where $B$ acts on $X$ via $\ominus J_n$) and the mapping spectrum formulation $\Map_{\Sp}(MTB, X)$. The equivalence between these formulations is a manifestation of the duality between tangential and normal framings, and we explain how the $\ominus$ sign in the GMTW formulation is absorbed by this duality. While the mapping spectrum perspective is well established in the literature \cite{MR2506750, MR4821355}, the homotopy fixed point perspective--and its direct connection to the $J$-homomorphism--appears to be less widely known, and we develop it here.

\subsection*{Outline}

This paper is organised as follows.
Section~\ref{sec:setup} recalls the background on bordism categories, the cobordism hypothesis, the $J$-homomorphism, and the localisation functor $L$.
Section~\ref{sec:localisation} constructs the $O(n)$-action on $L(\Bord_n^{\fr}) \simeq S$ via the abstract Pontryagin--Thom collapse map, proves that it is given by $\ominus J_n$ (Theorem~\ref{thm:main}), and identifies $L(\Bord_n^B)$ with the Madsen--Tillmann spectrum $MTB$ (Proposition~\ref{prop:MTB}).
Section~\ref{sec:computations} computes the induced action on the low-degree homotopy groups of $S$.
Section~\ref{sec:physics} discusses the physical interpretations in terms of invertible field theories, anomalies, and symmetry-protected topological phases.
\subsection*{Notation and conventions}

We introduce the notation to be used throughout the paper.

\begin{enumerate}
    \item[$\flat_1$] \emph{Framings.} A framing of a rank-$n$ vector bundle $V \to X$ is a trivialisation $V \cong \mathbb{R}^n$, equivalently a bundle isomorphism $\mathbb{R}^n \times X \to V$. A stable framing is a trivialisation of $V \oplus \mathbb{R}^N \cong \mathbb{R}^{n+N}$ for some $N \geq 0$. For a $k$-manifold $W$ with $n \geq k$, an $n$-framing is a trivialisation $TW \oplus \mathbb{R}^{n-k} \cong \mathbb{R}^n$.
    \item[$\flat_2$]\label{conventions-item:frame} \emph{Principal bundles and frame bundles.} Principal $G$-bundles have right actions. The frame bundle of a rank-$n$ vector bundle $V \to X$ is the principal $GL_n(\mathbb{R})$-bundle whose fibre over $x \in X$ is the collection of linear isomorphisms $T \colon \mathbb{R}^n \to V_x$, on which $A \in GL_n(\mathbb{R})$ acts on the right by precomposition: $T \cdot A := T \circ A$. A framing of $V$ is a section of the frame bundle, i.e.\ a map $\phi \colon\mathbb{R}^n \to V$ (the inverse of a trivialisation $V \to \mathbb{R}^n$). The right action of $g \in O(n)$ on framings is by precomposition: $g \cdot \phi := \phi \circ g$.

    \item[$\flat_3$] \emph{Madsen--Tillmann spectra.} For a tangential structure $\xi \colon B \to BO(n)$, the localisation of the bordism category is $L(\mathrm{Bord}_n^B) \simeq MTB$, where $MTB$ is the Madsen--Tillmann spectrum (Definition~\ref{def:MTB}).
   \item[$\flat_4$]\label{conventions-item:O-op} \emph{The identification $O(n)^{\mathrm{op}} \cong O(n)$.} For any $\infty$-group $G$, inversion $g \mapsto g^{-1}$ gives a canonical isomorphism $G \cong G^{\mathrm{op}}$. If $G$ is abelian, the identity also gives an isomorphism $G \cong G^{\mathrm{op}}$, and the two differ by the inversion automorphism. For $O(n)$ (non-abelian for $n \geq 2$), only the inversion isomorphism is available. For $O = \varinjlim_n O(n)$ (homotopy abelian), both are available, and they differ by an automorphism that we denote by $\ominus$. Throughout, we use the \emph{inversion} isomorphism $g \mapsto g^{-1}$ to identify $O(n)^{\mathrm{op}} \cong O(n)$ and $O^{\mathrm{op}} \cong O$, and this identification is compatible with the stabilisation maps $O(n) \to O(n+1)$ in the colimit. The commutativity identification is not used.
\end{enumerate}

\subsection*{Acknowledgements}
ASD and LS would like to thank Theo Johnson Freyd for reading a draft of the paper. ASD would like to thank Lennart Meier for enouraging him to think about problems arising in Physics, broadly speaking. ASD would also like to thank his PhD advisor Drew Heard and Shachar Carmeli for the support and encouragement. LS is grateful to Theo Johnson-Freyd and Lukas M\"uller for many discussions surrounding bordism $(\infty,n)$-categories and $O(n)$-actions.
\subsection*{LLM Declaration}
The authors did not use any artificial intelligence tools in the preparation of this manuscript.
\section{Setup}\label{sec:setup}
\noindent
In this section, we recall the essential objects and theorems that form the main background for our work: cobordism categories and the cobordism hypothesis, the $J$-homomorphism and Thom spectra, and the localisation functor from symmetric monoidal $(\infty,n)$-categories to spectra. Throughout, we work with the theory of $\infty$-categories, letting $\Cat_{(\infty,n)}$ denote the $\infty$-category of $(\infty,n)$-categories, and $\Sp$ the $\infty$-category of spectra. We use $\Fun^{\otimes}(-,-)$ to denote the space of symmetric monoidal functors. 

\subsection{Bordism categories and the cobordism hypothesis}\label{sec:bordism}

\begin{definition}
Let $\xi\colon B \to BO(n)$ be a map of spaces, which we call a \emph{tangential structure}. The \emph{bordism category} $\Bord_n^B$ is the symmetric monoidal $(\infty,n)$-category sketched as follows. For a complete definition, the reader is encouraged to consult \cite{MR3924174} and \cite[Section~2.2]{MR2555928}.

\begin{itemize}
    \item[$\natural_1$] \emph{Objects:} compact $0$-manifolds $P$ equipped with a $B$-structure on $TP \oplus \underline{\mathbb{R}}^n$. Since $TP$ is trivial (a $0$-manifold is a finite set of points), a $B$-structure is a map from each point to the homotopy fibre of $\xi$, i.e., a lift of the constant map to $\underline{\mathbb{R}}^n$ through $\xi$.
    
    \item[$\natural_2$] \emph{$1$-morphisms:} a $1$-morphism from $P_0$ to $P_1$ is a compact $1$-manifold $W$ with a $B$-structure on $TW \oplus \underline{\mathbb{R}}^{n-1}$ and an identification $\partial W = P_0 \sqcup P_1$ compatible with the $B$-structures using in- and outgoing normal vectors on $P_0$ and $P_1$ respectively.
    
    \item[$\natural_3$] \emph{$k$-morphisms, $2 \le k \le n$:} a $k$-morphism between two $(k-1)$-morphisms is a compact $k$-manifold $W$ with corners, a $B$-structure on $TW \oplus \underline{\mathbb{R}}^{n-k}$ and an identification of its boundary with the incoming and outgoing $(k-1)$-morphisms.
    
    \item[$\natural_4$] \emph{$(n+1)$-morphisms:} diffeomorphisms of $n$-manifolds preserving the $B$-structure.
    \item[$\natural_5$] \emph{Higher morphisms:} Isotopies of $B$-diffeomorphisms.
    
    \item[$\natural_6$] \emph{Composition:} glueing of bordisms along boundaries.
\end{itemize}
The symmetric monoidal structure is given by disjoint union.
\end{definition}

\begin{definition}
For $n \geq k$, an $n$-framing on a $k$-manifold $W$ (potentially with corners) is a $B$-structure on $TW \oplus \underline{\mathbb{R}}^{n-k}$ for $B = \pt$.
\end{definition}

Recall that the \emph{frame bundle} of an $n$-dimensional vector bundle $V \to X$ is the principal $GL_n(\R)$-bundle whose fibre over $x \in X$ is the collection of linear isomorphisms $T \colon \R^n \to V_x$ on which $A \in GL_n(\R)$ acts by right composition. 
A section of the frame bundle is a basis of all the $V_x$ continuous in $x$.
An $n$-framing on $W$ is a trivialisation of $TW \oplus \underline{\mathbb{R}}^{n-k}$, which therefore corresponds to a trivialisation of its frame bundle.

\begin{remark}
    The data of an $n$-framing on a $k$-manifold for $n>k$ is typically weaker than a $k$-framing.
    For example, $S^k$ always has a $(k+1)$-framing but only a $k$-framing for $k = 0,1,3$ and $7$.
\end{remark}

\begin{definition}
When $B = \pt$ with the unique map $\pt \to BO(n)$, we obtain the \emph{fully framed bordism category} $\Bord_n^{\fr}$. 
\end{definition}

\begin{remark}
\label{rmk:O(n)-action}
The $GL_n(\mathbb{R})$-action on the frame bundle yields a continuous right action of $O(n)$ on $\mathrm{Bord}_n^{\mathrm{fr}}$ by change of framing. More precisely, there is a map of $\mathbb{E}_1$-spaces
\[
\alpha_{\mathrm{CH}} \colon O(n)^{\mathrm{op}} \to \mathrm{Aut}_{\mathrm{Cat}_{(\infty,n)}^{\mathbb{E}_\infty}}(\mathrm{Bord}_n^{\mathrm{fr}}).
\]
Alternatively, one can describe the delooping as a map
\[
BO(n) \to B\mathrm{Aut}(\mathrm{Bord}_n^{\mathrm{fr}}) \hookrightarrow \mathrm{Cat}_{(\infty,n)}^{\mathbb{E}_\infty},
\]
which is the input for the cobordism hypothesis.
\end{remark}

\noindent
A \emph{fully extended framed $n$-dimensional topological field theory} valued in a symmetric monoidal $(\infty,n)$-category $\mathcal{C}$ is a symmetric monoidal functor $Z\colon\Bord_n^{\fr}\to \mathcal{C}$. We denote the space of such functors by $\TFT^{\fr}_n(\mathcal{C}) := \Fun^{\otimes}(\Bord_n^{\fr}, \mathcal{C})$. The foundational result in this subject is the cobordism hypothesis, stated by Baez-Dolan \cite{baez1995higher}.

\begin{theorem}[Cobordism Hypothesis, \cite{MR2555928, grady2021geometric}]\label{thm:cobordism-hypothesis}
Let $\mathcal{C}$ be a symmetric monoidal $(\infty, n)$-category with duals. Evaluation at the positively framed point $\pt^+ \in \Bord_n^{\fr}$ induces an equivalence of spaces
\[
\TFT_n^{\fr}(\mathcal{C}) \xrightarrow{\;\sim\;} \mathcal{C}^{\fd},
\]
where $\mathcal{C}^{\fd}$ denotes the underlying $\infty$-groupoid of fully dualisable objects in $\mathcal{C}$. 
\end{theorem}

\begin{remark}\label{rem:CH-action}
The above theorem gives an $O(n)$-action on $\mathcal{C}^{\fd}$ from the $O(n)$-action on $\Bord_n^{\fr}$ via the equivalence.
The $\infty$-category $\Fun(BO(n),\mathcal{S})$ of spaces with an $O(n)$-action remembers all higher coherence of the $O(n)$-equivariance. 
Concretely, the equivalence $\TFT_n^{\fr}(\mathcal{C}) \simeq \mathcal{C}^{\fd}$ interwines the two $O(n)$-actions via a specified homotopy, together with a compatible family of higher homotopies witnessing the $O(n)$-coherence. At the level of points, this implies that for any $g \in O(n)$, the diagram \[\begin{tikzcd}
    \TFT^{\fr}_n(\mathcal{C}) \ar[r, "\sim"] \ar[d, "g\cdot"'] & \mathcal{C}^{\fd} \ar[d, "g\cdot"]\\
    \TFT_n^{\fr}(\mathcal{C}) \ar[r, "\sim"] & \mathcal{C}^{\fd} 
\end{tikzcd}\] commutes, where the left vertical map precomposes with the change-of-framing functor $\Bord_n^{\fr}\xrightarrow[]{g}\Bord^{\fr}_n$, and the right vertical map is the action constructed by Lurie \cite{MR2555928}
\end{remark}

Using the equivalence $\Bord_n^B \simeq \colim_{B} \Bord_n^{\fr}$ (the colimit of the $B$-action via $\xi$) and the fact that $\Fun^{\otimes}(-,\mathcal{C})$ converts colimits to limits, one gets the corollary:
\begin{corollary}\label{cor:CH-tangential}
For any $\mathcal{C}$ as above and any map $\xi\colon B \to BO(n)$,
\[
\TFT_n^B(\mathcal{C}) \simeq (\mathcal{C}^{\fd})^{hB},
\]
where $(-)^{hB}$ denotes homotopy fixed points for the action of $B$ on $\mathcal{C}^{\fd}$ through the composition $B \xrightarrow{\xi} BO(n) \to \Aut(\mathcal{C}^{\fd})$.
\end{corollary}

\subsection{Units, the \texorpdfstring{$J$}{J}-homomorphism, and Thom spectra}\label{sec:units-J}

We denote by $\Sp$ the symmetric monoidal $\infty$-category of spectra. Its monoidal unit is the sphere spectrum $S$. By a classical theorem of May \cite[Theorem 14.4, Corollary 14.5]{may2007geometry}, and its $\infty$-categorical refinement \cite[Remark 5.2.6.26]{Lurie2017HA}, there is an equivalence between the $\infty$-category of connective spectra $\Sp_{\ge 0}$ and the $\infty$-category of grouplike $\E_\infty$-spaces. For a spectrum $X$, we let $\Omega^\infty X \in \mathcal{S}$ denote its underlying (infinite loop) space.
Recall that an $\E_\infty$-ring spectrum is a commutative algebra object in $\Sp$.

\begin{definition}
\label{def:units}
Let $R$ be an $\E_\infty$-ring spectrum.
\begin{itemize}
    \item[$\natural_1$] the \emph{higher group of units} $\gl_1(R)$ is defined as the $\E_\infty$-group of automorphisms of the unit $R \in \Mod_R$ in the category of right $R$-modules.
    \item[$\natural_2$] The \emph{Picard groupoid $\Pic(R)$ of $R$} is the $\E_\infty$-group of invertible modules and invertible module maps between them.
\end{itemize}
\end{definition}

We will typically identify $\gl_1(R)$ and $\Pic(R)$ with their corresponding connective spectra.
We have $\Omega \Pic(R) = \gl_1(R)$.
If $\GL_1(R) := \Omega^\infty \gl_1(R)$ is the space of units~\cite[Definition 2.1]{MR3286898}, there is pullback
\[
\begin{tikzcd}
\GL_1(R) \arrow[r] \arrow[d] & \Omega^\infty R \arrow[d, "\pi_0"] \\
(\pi_0 R)^\times \arrow[r] & \pi_0 R
\end{tikzcd}
\]
in the $\infty$-category $\mathcal{S}$~\cite[Proposition 6.2]{ando2009unitsringspectrathom}.

Since $S$ is the monoidal unit in $\Sp$, it has a canonical $\E_\infty$-structure and $gl_1(S) = \Aut_{\Sp}(S)$.
For the sphere spectrum $S$, we have $\pi_0(S) \cong \mathbb{Z}$ and therefore $(\pi_0 S)^\times \cong \{\pm 1\} \cong \mathbb{Z}/2$. The space $\GL_1(S)$ has two connected components, and its homotopy groups are
\[
\pi_k(\GL_1(S)) \cong
\begin{cases}
\mathbb{Z}/2 & \text{for } k = 0, \\
\pi_k(S) & \text{for } k \ge 1.
\end{cases}
\]

\begin{definition}
\label{def:Jn}
For a fixed dimension $n$, the \emph{unstable $J$-homomorphism} is the map of $\mathbb{E}_1$-spaces
\[
J_n \colon O(n) \longrightarrow \mathrm{GL}_1(S)
\]
defined as follows. A point $g \in O(n)$ defines a linear isometry $\mathbb{R}^n \to \mathbb{R}^n$, whose one-point compactification is a based map $S^n \to S^n$, an element of the $n$-fold loop space $\Omega^n S^n$. We set $J_n(g)$ to be the image under
\[
O(n) \hookrightarrow \Omega^n S^n \xrightarrow{\ \text{stab}\ } \Omega^\infty S,
\]
where the first map is the classical unstable $J$-homomorphism \cite[pp.~291--293]{Atiyah1961} and $\text{stab}$ is the stabilisation map $\Omega^n S^n \to \varinjlim_k \Omega^k S^k =: \Omega^\infty S$. The image of the composition lies in the units $\mathrm{GL}_1(S) \subset \Omega^\infty S$ by inspection of $\pi_0$: an element $g \in O(n)$ maps to a based map $S^n \to S^n$ of degree $\det(g) = \pm 1$, whose stabilisation has degree $\pm 1$ in $\pi_0(S) \cong \mathbb{Z}$, hence lies in the units $\pi_0(S)^\times = \{\pm 1\} \subset \mathrm{GL}_1(S)$. Note that the stabilisation map $\Omega^n S^n \to \Omega^\infty S$ is \emph{not} a map to $\mathrm{GL}_1(S)$ on all of $\Omega^n S^n$: the components of degree $\neq \pm 1$ do not land in the units.
\end{definition}

The unstable $J$-homomorphisms $J_n$ are compatible with the stabilisation inclusions $O(n) \hookrightarrow O(n+1)$, in the sense that the diagram
\[
\begin{tikzcd}
O(n) \arrow[r, "J_n"] \arrow[d, hook] & \mathrm{GL}_1(S) \arrow[d, equal] \\
O(n+1) \arrow[r, "J_{n+1}"] & \mathrm{GL}_1(S)
\end{tikzcd}
\]
commutes up to homotopy. Passing to the colimit, we obtain the \emph{stable $J$-homomorphism}
\[
J \colon O = \varinjlim_n O(n) \longrightarrow \mathrm{GL}_1(S).
\]
Equivalently, at the level of classifying spaces, this gives $BJ \colon BO \to B\mathrm{GL}_1(S)$, which classifies stable spherical fibrations.

For the purpose of defining Thom spectra, it is convenient to give a more conceptual construction of $J$ that packages the whole family $\{J_n\}_n$ into a single map. Let $\mathrm{Vec} \simeq \bigsqcup_{n \geq 0} BO(n)$ denote the $\mathbb{E}_\infty$-groupoid of real vector spaces and orthogonal maps, under Whitney sum $\oplus$. One-point compactification with basepoint at infinity defines a symmetric monoidal functor
\begin{equation}
\label{eq:one-point-compactification}
\mathrm{Vec} \to \mathcal{S}_*, \qquad V \mapsto S^V,
\end{equation}
where we give pointed spaces the Cartesian symmetric monoidal structure. Composing \eqref{eq:one-point-compactification} with the symmetric monoidal left adjoint $\Sigma^\infty : \mathcal{S}_* \to \mathrm{Sp}$ of $\Omega^\infty$, we obtain a symmetric monoidal functor
\[
(\mathrm{Vec}, \oplus) \to (\mathrm{Sp}, \wedge).
\]
Since $S^V \wedge S^{-V} \simeq S$ and $O(V)$ maps to automorphisms of $S^V$, this functor lands in the non-full subcategory $\mathrm{Pic}(S) \hookrightarrow \mathrm{Sp}$ of invertible $S$-modules. Since $\mathrm{Pic}(S)$ is grouplike, the universal property of the group completion $ko = BO \times \mathbb{Z}$ of $\mathrm{Vec}$ implies that the composition factors through a map
\[
J \colon BO \times \mathbb{Z} \longrightarrow \mathrm{Pic}(S).
\]
The restriction of $J$ to the component $BO(n) \subset BO \times \mathbb{Z}$ recovers the unstable $J_n$: on the level of $\pi_0$, both send a rank-$n$ vector space to $\Sigma^n S$, and the induced action of $O(n)$ on $\Sigma^n S$ is the same. Thus the conceptual construction above is a repackaging of the family $\{J_n\}_n$ into a single stable map.

\begin{definition}
\label{def:stable-J}
The \emph{stable $J$-homomorphism} is the map
\[
J \colon BO \times \mathbb{Z} \longrightarrow \mathrm{Pic}(S)
\]
constructed above.
\end{definition}

\begin{definition}
\label{def:stable-vector-bundle}
A \emph{stable vector bundle} over $X$ is a map $X \to BO \times \mathbb{Z}$. If $X \to BO(n)$ is a map classifying a rank-$n$ vector bundle, its \emph{stabilisation} is its composition with the group completion map $BO(n) \to BO \times \mathbb{Z}$ that includes into the $\{n\} \subseteq \mathbb{Z}$ component.
\end{definition}

\begin{remark}
\label{rmk:J-on-virtual-bundles}
Observe that in our definition, a rank-$n$ vector bundle includes into the $\{n\} \subseteq \mathbb{Z}$ component, and stable vector bundles can have a nonzero virtual rank. Because $J$ is defined using the universal property of the group completion, we have $J(V - W) = S^V \wedge S^{-W}$ for vector spaces $V, W$.
\end{remark}

\begin{definition}
\label{def:Thom}
The \emph{Thom spectrum} $X^V$ of a stable vector bundle $X \to BO \times \mathbb{Z}$ is the colimit of
\[
X \to BO \times \mathbb{Z} \xrightarrow{\ J\ } \mathrm{Pic}(S) \to \mathrm{Sp}.
\]
\end{definition}

\begin{remark}
\label{rmk:Thom-as-cofib}
A functor $X \to \mathcal{S}_*$ corresponds to a fibration with total space $E$, its colimit. For $X \to \mathrm{Vec}$ a vector bundle, the Thom space can be computed as
\[
\operatorname*{colim}(X \to \mathcal{S}_* \to \mathcal{S}) \sqcup_X \mathrm{pt} = \operatorname*{cofib}(s \colon X \to E),
\]
where $E = \operatorname*{colim}(X \to \mathcal{S}_* \to \mathcal{S})$ is the total space of the associated spherical fibration, and $s$ its canonical section. This connects Definition~\ref{def:Thom} with the usual definition of Thom space.
\end{remark}

\begin{definition}
\label{def:MTB}
For a tangential structure $\xi \colon B \to BO(n)$, the \emph{Madsen--Tillmann spectrum} $MTB$ is the Thom spectrum of the virtual normal bundle, i.e.\ the colimit of
\[
B \xrightarrow{\ \ominus(\iota, n) \circ \xi\ } BO \times \mathbb{Z} \xrightarrow{\ J\ } \mathrm{Pic}(S) \to \mathrm{Sp},
\]
where $\iota \colon BO(n) \to BO$ is the stabilisation map, $n \colon BO(n) \to \mathbb{Z}$ is the constant map at $n$, and $\ominus$ is the inverse in the grouplike $\mathbb{E}_\infty$-structure on $BO \times \mathbb{Z}$ coming from Whitney sum. Equivalently, under the identification of Definition~\ref{def:Thom}, $MTB = B^{-\xi}$, the Thom spectrum of the virtual bundle $-\xi$ of rank $-n$. A closed $n$-manifold $M$ with a $B$-structure gives a map $S \to MTB$ via the Pontryagin--Thom construction, so that $MTB$ is the spectrum of $B$-structured bordism classes in degree $0$.
\end{definition}
\begin{example}
\label{ex:trivial-bundle}
Note that $BO(n) \to \mathrm{Pic}(S)$ maps to $\Sigma^n S$. Therefore, the Thom spectrum of an $n$-dimensional trivial bundle over $X$ is the colimit over the trivial diagram indexed by $X$, which is $\Sigma^n \Sigma^\infty_+ X$.
\end{example}
\subsection{The localisation functor \texorpdfstring{$L$}{L}}\label{sec:localisation-functor}
\subsubsection{The $\infty$-category of $(\infty,n)$-categories and the localisation functor}\label{sec:Cat-infinity-n}

We recall the recursive definition of the $\infty$-category of $(\infty,n)$-categories and the adjunctions between the various levels. The primary reference is \cite{GH19}.

\begin{definition}\label{def:Cat-infinity-n}
We define $\mathrm{Cat}_{(\infty,n)}$ recursively by
\[
\mathrm{Cat}_{(\infty,0)} := \mathcal{S},\qquad \mathrm{Cat}_{(\infty,n)}:=\mathrm{Cat}[\mathrm{Cat}_{(\infty,n-1)}],
\]
where $\mathcal{S}$ is the $\infty$-category of spaces and $\mathrm{Cat}[-]$ is the construction of \cite[\S 5]{GH19}  taking a monoidal $\infty$-category $\mathcal{V}$ to the $\infty$-category $\mathrm{Cat}[\mathcal{V}]$ of small $\mathcal{V}$-enriched $\infty$-categories.
    
\end{definition}

By \cite[Theorem 5.4.6]{GH19}, we have $\mathrm{Cat}_{(\infty,1)} \simeq \mathrm{Cat}_\infty$, and the construction is Cartesian presentably symmetric monoidal at every level \cite[Proposition 5.7.16]{GH19}.

There are two natural adjunctions between consecutive levels:
\begin{enumerate}
\item[$\natural_1$] \textit{Core / inclusion.} The \emph{core functor}
\[
(-)^{\mathrm{core}_{n-1}} \colon \mathrm{Cat}_{(\infty,n)} \to \mathrm{Cat}_{(\infty,n-1)}
\]
is right adjoint to the fully faithful inclusion
\[
\iota_n \colon \mathrm{Cat}_{(\infty,n-1)} \hookrightarrow \mathrm{Cat}_{(\infty,n)}
\]
of $(\infty,n-1)$-categories as $(\infty,n)$-categories with only invertible $n$-morphisms.
\item[$\natural_2$] \textit{Localisation / inclusion.} The \emph{groupoid-completion functor}
\[
L_{n-1} \colon \mathrm{Cat}_{(\infty,n)} \to \mathrm{Cat}_{(\infty,n-1)}
\]
is left adjoint to $\iota_n$. It universally inverts all $n$-morphisms.
\end{enumerate}

The composite \(L_0 \colon \mathrm{Cat}_{(\infty,n)} \to \mathcal{S}\) is the $\infty$-groupoid completion functor. 

\begin{lemma}
\label{lem:Cat-preserves-adjunctions}
The construction $\mathrm{Cat}[-]$ preserves adjunctions: if $F \colon \mathcal{C} \rightleftarrows \mathcal{D} \colon G$ is an adjunction of monoidal $\infty$-categories, then
\[
\mathrm{Cat}[F] \colon \mathrm{Cat}[\mathcal{C}] \rightleftarrows \mathrm{Cat}[\mathcal{D}] \colon \mathrm{Cat}[G]
\]
is an adjunction of $\infty$-categories.
\end{lemma}

\begin{proof}
By \cite[Proposition 5.7.8 and Proposition 5.7.16]{GH19}, the construction $\mathrm{Cat}[-]$ assembles into a functor
\[
\mathrm{Cat}[-] \colon \mathrm{CAlg}(\mathrm{Pr}^{\mathrm{L}}) \to \mathrm{CAlg}(\mathrm{Pr}^{\mathrm{L}}).
\]
Since $\mathrm{Pr}^{\mathrm{L}}$ is the $\infty$-category of presentable $\infty$-categories and left adjoints, such a functor sends left adjoints to left adjoints and right adjoints to right adjoints. Therefore, it preserves adjunctions.
\end{proof}

\begin{lemma}
\label{lem:Cat-preserves-products}
The construction $\mathrm{Cat}[-]$ preserves finite products: if $\mathcal{C}, \mathcal{D}$ are monoidal $\infty$-categories, then
\[
\mathrm{Cat}[\mathcal{C} \times \mathcal{D}] \simeq \mathrm{Cat}[\mathcal{C}] \times \mathrm{Cat}[\mathcal{D}].
\]
More generally, if $\mathcal{C} \to \mathcal{D}$ is a symmetric monoidal functor, then so is $\mathrm{Cat}[\mathcal{C}] \to \mathrm{Cat}[\mathcal{D}]$.
\end{lemma}

\begin{proof}
By \cite[Proposition 5.7.16]{GH19}, the construction $\mathrm{Cat}[-]$ is symmetric monoidal with respect to the Cartesian symmetric monoidal structure on $\mathrm{Cat}_{(\infty,n-1)}$. Since the monoidal structure is Cartesian, preserving it means preserving products. The second statement follows from the same citation.
\end{proof}

\begin{remark}
\label{rmk:inductive-sym-mon}
Combining Lemmas~\ref{lem:Cat-preserves-adjunctions} and \ref{lem:Cat-preserves-products}, we obtain by induction that:
\begin{enumerate}
\item[$\natural_1$] The adjunctions $L_{n-1} \dashv \iota_n$ are symmetric monoidal adjunctions, for every $n \geq 1$.
\item[$\natural_2$] The $\infty$-groupoid completion functor $L_0 \colon \mathrm{Cat}_{(\infty,n)} \to \mathcal{S}$ is symmetric monoidal.
\end{enumerate}
This is the key inductive input for the symmetric monoidality of the localisation functor $L$ (Lemma~\ref{lem:L-symmetric-monoidal}).
\end{remark}

Let $\Cat_{(\infty,n)}^{\E_\infty}$ denote the $\infty$-category of symmetric monoidal $(\infty,n)$-categories and symmetric monoidal functors between them. There is a natural fully faithful functor
\[
\iota\colon \Sp_{\ge 0} \hookrightarrow \Cat_{(\infty,n)}^{\E_\infty}
\]
defined as follows. A connective spectrum $X$ determines a grouplike $\E_\infty$-space $\Omega^\infty X$, which is in particular a symmetric monoidal $(\infty,0)$-category (an $\infty$-groupoid). We view this as a symmetric monoidal $(\infty,n)$-category with only invertible $k$-morphisms for all $1 \le k \le n$. Under the homotopy hypothesis, this is the $\infty$-category obtained from $\Omega^\infty X$ by taking the fundamental $\infty$-groupoid in which paths as $1$-morphisms, homotopies as $2$-morphisms, and so on. All objects of $\iota(X)$ are fully dualisable because the monoidal structure is grouplike: every object $x$ has an inverse $x^{-1}$ under the tensor product, which serves as its dual. 

The $\infty$-category $\Cat_{(\infty,n)}^{\E_\infty}$ is presentable \cite[Proposition~5.1.7 + Definition~5.3.6]{stefanich2020presentable}, and $\iota$ is an accessible right adjoint. Therefore, by the adjoint functor theorem for presentable $\infty$-categories \cite[Theorem~4.1.1]{nguyen2020adjoint}, $\iota$ admits a left adjoint
\[
L\colon \Cat_{(\infty,n)}^{\E_\infty} \longrightarrow \Sp_{\ge 0}.
\]
We call $L$ the \emph{localisation functor}. On a symmetric monoidal $(\infty,n)$-category $\mathcal{C}$, the functor $L$ universally forces all $k$-morphisms for $0 \le k \le n$ to become equivalences.
Since $k=0$ is included, this also performs group completion with respect to the symmetric monoidal structure. 
Concretely, the underlying infinite loop space of $L(\mathcal{C})$ is the group completion of the localisation $L_0 \mathcal{C}$ of $\mathcal{C}$ with respect to $k$-morphisms for $0 < k \leq n$. 
We write
\[
L_k \colon \mathrm{Cat}_{(\infty,n)}^{\E_\infty}  \longrightarrow \mathrm{Cat}_{(\infty,k)}^{\E_\infty} 
\]
for the level-wise localisation functor (for $k < n$), i.e.\ the functor that inverts all $j$-morphisms for $j > k$, and
\[
L \colon \mathrm{Cat}_{(\infty,n)}^{\mathbb{E}_\infty} \longrightarrow \mathrm{Sp}_{\geq 0}
\]
for the global localisation functor, left adjoint to $\iota$ as above. The global $L$ factors as the composite
\[
L \colon \mathrm{Cat}_{(\infty,n)}^{\mathbb{E}_\infty} \xrightarrow{\ L_{n-1}\ } \mathrm{Cat}_{(\infty,n-1)}^{\mathbb{E}_\infty} \xrightarrow{\ L_{n-2}\ } \cdots \xrightarrow{\ L_0\ } \mathrm{Cat}_{(\infty,0)}^{\mathbb{E}_\infty} = \mathcal{S}^{\mathbb{E}_\infty} \xrightarrow{\ \text{group completion}\ } \mathrm{Sp}_{\geq 0},
\]
where the intermediate localisations are the level-wise ones restricted to the symmetric monoidal subcategory.

\begin{lemma}
\label{lem:Lk-preserves-products}
Let $k < n$. The level-wise localisation $L_k \colon \mathrm{Cat}_{(\infty,n)} \to \mathrm{Cat}_{(\infty,k)}$ preserves finite products.
\end{lemma}

\begin{proof}
The $\infty$-category $\mathrm{Cat}_{(\infty,n)}$ has exponential objects given by functor $(\infty,n)$-categories: if $\mathcal{C} \in \mathrm{Cat}_{(\infty,n)}$ and $\mathcal{D} \in \mathrm{Cat}_{(\infty,k)}$, then $\mathcal{C}^{\mathcal{D}} \in \mathrm{Cat}_{(\infty,k)}$ (the exponentials are computed levelwise, and $L_k$ only affects the top levels). It follows for $\mathcal{C}_1, \mathcal{C}_2 \in \mathrm{Cat}_{(\infty,n)}$ and $\mathcal{D} \in \mathrm{Cat}_{(\infty,k)}$ that
\[
\mathrm{Hom}_{\mathrm{Cat}_{(\infty,k)}}(L_k(\mathcal{C}_1 \times \mathcal{C}_2), \mathcal{D}) \cong \mathrm{Hom}_{\mathrm{Cat}_{(\infty,k)}}(L_k \mathcal{C}_1, \mathcal{D}^{\mathcal{C}_2}) \cong \mathrm{Hom}_{\mathrm{Cat}_{(\infty,n)}}(\mathcal{C}_1, \mathcal{D}^{\mathcal{C}_2}) \cong \mathrm{Hom}_{\mathrm{Cat}_{(\infty,n)}}(\mathcal{C}_1 \times \mathcal{C}_2, \mathcal{D}).
\]
By Yoneda, the unit map $\mathcal{C}_1 \times \mathcal{C}_2 \to L_k \mathcal{C}_1 \times \mathcal{C}_2$ becomes an equivalence after applying $L_k$. Repeating the argument for the other factor proves the same for $\mathcal{C}_1 \times \mathcal{C}_2 \to L_k \mathcal{C}_1 \times L_k \mathcal{C}_2$. Since $L_k \mathcal{C}_1 \times L_k \mathcal{C}_2 \in \mathrm{Cat}_{(\infty,k)}$, we obtain
\[
L_k(\mathcal{C}_1 \times \mathcal{C}_2) \cong L_k(L_k \mathcal{C}_1 \times L_k \mathcal{C}_2) \cong L_k \mathcal{C}_1 \times L_k \mathcal{C}_2,
\]
so $L_k$ preserves binary products. The terminal object is preserved since $L_k(\mathrm{pt}) \cong \mathrm{pt}$ (the terminal $(\infty,n)$-category has only one object and only identity morphisms, hence is already local). Therefore $L_k$ preserves finite products.
\end{proof}

\begin{lemma}
\label{lem:L-symmetric-monoidal}
The global localisation functor $L \colon \mathrm{Cat}_{(\infty,n)}^{\mathbb{E}_\infty} \to \mathrm{Sp}_{\geq 0}$ is symmetric monoidal, with respect to the Cartesian symmetric monoidal structure on $\mathrm{Cat}_{(\infty,n)}^{\mathbb{E}_\infty}$ and the smash product on $\mathrm{Sp}_{\geq 0}$.
\end{lemma}

\begin{proof}
The functor $\iota \colon \mathrm{Sp}_{\geq 0} \hookrightarrow \mathrm{Cat}_{(\infty,n)}^{\mathbb{E}_\infty}$ factors as
\[
\mathrm{Sp}_{\geq 0} \xrightarrow{\ \iota_1\ } \mathcal{S} \xrightarrow{\ \iota_2\ } \mathrm{Cat}_{(\infty,n)}^{\mathbb{E}_\infty},
\]
where $\iota_1$ includes connective spectra as grouplike $\mathbb{E}_\infty$-spaces and $\iota_2$ includes spaces as $(\infty,n)$-categories with only invertible morphisms. The left adjoint to $\iota_2$ is the $\infty$-groupoid completion functor, which is symmetric monoidal by Remark~\ref{rmk:inductive-sym-mon}. The left adjoint to $\iota_1$ is group completion, which is symmetric monoidal (group completion is a colimit of symmetric monoidal functors, so it inherits a symmetric monoidal structure). Therefore, $L$, being the composite of these two symmetric monoidal left adjoints, is symmetric monoidal.
\end{proof}

\begin{remark}\label{rem:L-duals}
If $\mathcal{C}$ is a symmetric monoidal $(\infty,n)$-category with duals, then for every $k > 0$, any $k$-morphism has a left and a right adjoint \cite[Definition~2.3.13+Definition~2.3.16]{MR2555928} and is therefore already invertible in $L_1(\mathcal{C})$. Consequently, $L_1(\mathcal{C})$  is already the full localisation $L(\mathcal{C})$. This is precisely the situation for the bordism category $\Bord_n^{B}$ when $n > 0$.
\end{remark}
\section{Localisation of the Cobordism Hypothesis Action}\label{sec:localisation}

\subsection{Recovering the sphere spectrum from abstract nonsense}
\label{sec:abstract}

One way to phrase the cobordism hypothesis is to say that $\Bord^{\fr}_n$ is the free symmetric monoidal $(\infty,n)$-category with duals on a single object.
Geometrically, this object is the connected $0$-dimensional manifold with its positive framing.
Since the localisation functor $L$ is a left adjoint, it preserves colimits, and therefore preserves the property of being free.
Therefore, we expect $L(\Bord^{\fr}_n)$ to be the free connective spectrum on a point, i.e.\ the sphere spectrum:

\begin{proposition}\label{prop:L-S}
    There is a unique equivalence of spectra $L(\Bord^{\fr}_n) \cong S$ with the property that the positively framed point gets mapped to $1 \in \pi_0 [S,S] \cong \pi_0(S)$.
\end{proposition}
\begin{proof} 
For any connective spectrum $X$, we compute:
\begin{equation}
\label{eq:invcobhyp}
\Hom_{\Sp_{\ge 0}}(L(\Bord_n^{\fr}), X) \xrightarrow{(1)} \Hom_{\Cat_{(\infty,n)}^{\E_\infty}}(\Bord_n^{\fr}, \iota(X)) \xrightarrow{(2)} \iota(X)^{\fd} \xrightarrow{(3)} \Omega^\infty X.
\end{equation}
Here:
\begin{enumerate}
    \item[$\natural_1$] is the $L\dashv \iota$ adjunction.
    \item[$\natural_2$] is the cobordism hypothesis (Theorem \ref{thm:cobordism-hypothesis}), which states that symmetric monoidal functors $\Bord_n^{\fr} \to \iota(X)$ correspond to fully dualisable objects of $\iota(X)$.
    \item[$\natural_3$] is the observation that $\iota(X)^{\fd}$, the $\infty$-groupoid of fully dualisable objects of $\iota(X)$ is equivalent to the underlying space $\Omega^\infty X$. Indeed, every object of $\iota(X)$ is fully dualisable, so the core is exactly the space of objects.
\end{enumerate}
The entire chain is natural in $X$. Thus $L(\Bord_n^{\fr})$ represents the $\infty$-groupoid-valued functor $X \mapsto \Omega^\infty X$ on $\Sp_{\ge 0}$. The sphere spectrum $S$ satisfies the same universal property, since $\Hom_{\Sp}(S, X) \simeq \Omega^\infty X$ for any spectrum $X$. By the $\infty$-categorical Yoneda lemma, we obtain an equivalence $L(\Bord_n^{\fr}) \simeq S$.

Explicitly, the equivalence is determined by the image of the generating object $x$: the unique map $L(\Bord_n^{\fr}) \to S$ sending $x$ to $1 \in \pi_0(S) \subset \Omega^\infty S$ is an equivalence.
Any other equivalence $L(\Bord_n^{\fr}) \cong S$ would differ by an automorphism of $S$.
Since $\pi_0 \Aut(S) = \pi_0 [S,S]^\times = \Z^\times = \{\pm 1\}$, this would change $L(\Bord_n^{\fr}) \to S$ by mapping the positively framed point to $-1 \in S$.
\end{proof}

\begin{remark}
    A similar argument appears in \cite[Section 7]{MR2994995}.
\end{remark}

\begin{remark}
Since $S$ is the monoidal unit in $\Sp_{\ge 0}$, it is the unique spectrum so that $\Hom_{\Sp_{\ge 0}}(S, X) \cong X$ for all $X$; any other such spectrum $Y$ would be right adjoint to $Y \wedge (-)$. 
We thus obtain the equivalence\[L(\Bord_n^{\fr})\simeq S\] of the above proposition also from uniqueness of adjoint functors.

As a consequence, the equivalence is given by 
\[
L(\Bord_n^{\fr}) \to L(\Bord_n^{\fr}) \wedge \Hom_{\Sp_{\geq 0}}(L(\Bord_n^{\fr}), S) \xrightarrow{ev} S
\]
where the first map is the inverse of \eqref{eq:invcobhyp} for $X = S$.
In other words, it evaluates a framed cobordism $M$ on the unique TQFT $\Bord_n \to S$ that assigns $1$ to the positively framed point.

The inverse map $S \to L(\Bord_n^{\fr})$ is given by the composition
\[
S \to S \wedge \Hom_{\Sp}(L(\Bord_n^{\fr}),L(\Bord_n^{\fr})) \to \Hom_{\Sp}(L(\Bord_n^{\fr}),L(\Bord_n^{\fr})) \to L(\Bord_n^{\fr}),
\]
where the first map is induced by the identity on $L(\Bord_n^{\fr})$, the second map is the canonical $S$-module structure that exists on any spectrum, and the final map is \eqref{eq:invcobhyp} applied to $X = L(\Bord_n^{\fr})$.
\end{remark}

\subsection{Tangential versus normal framings}
\label{sec:tangential-normal}

\begin{definition}
\label{def:framing}
Let $\xi \colon X \to BO(n)$ be a rank-$n$ vector bundle over $X$.
\begin{enumerate}
\item[$\natural_1$] An \emph{unstable framing} of $\xi$ is a nullhomotopy of $\xi$, i.e.\ an isomorphism of vector bundles $\mathbb{R}^n \times X \cong \xi$.
\item[$\natural_2$] A \emph{stable framing} of $\xi$ is a nullhomotopy of the composition $X \xrightarrow{\xi} BO(n) \to BO$.
\end{enumerate}
A stable (resp.\ unstable) framing of a manifold (with corners) $M$ is a stable (resp.\ unstable) framing of its tangent bundle $\tau_M$.
\end{definition}

If $X$ is compact, a stable framing is equivalent to a trivialisation of the bundle $\xi \oplus \mathbb{R}^N$ for $N$ sufficiently large. Since a stably trivial bundle of rank $> d$ over a $k$-dimensional CW-complex is trivial, it suffices to take $N = 1$ if $M$ is a manifold and $\xi = \tau_M$. However, the space of stable framings is a torsor over the $\infty$-group $\Omega \mathrm{Map}(M, BO) \cong \mathrm{Map}(M, O)$, which is in general different from $\mathrm{Map}(M, O(n+1))$.

If $j \colon M^n \hookrightarrow N^{n+k}$ is an embedding between manifolds, the $k$-dimensional quotient bundle $\nu_j := j^* \tau_N / \tau_M$ is called the \emph{normal bundle}. Since every exact sequence of vector bundles splits noncanonically, in case $N = \mathbb{R}^{n+k}$, this picks out a bundle $\nu_j$ such that $\nu_j \oplus \tau_M \cong \mathbb{R}^{n+k}$. By the Whitney embedding theorem, such $j$ always exists. Moreover, the induced stable normal bundle $M \to BO(k) \to BO$ is independent of the choice of $j$ and will be denoted $\nu_M$. This stable bundle is a choice of inverse of the stable tangent bundle under $\oplus$.

Since $\ominus \tau_M = \nu_M$, a stable framing is equivalent to a trivialisation of the stable normal bundle, but the $O$-actions differ by $\ominus$. More geometrically, suppose we have a given $j \colon M \to \mathbb{R}^{n+k}$ giving, after a contractible choice of metric, a trivialisation $\nu_j \oplus \tau_M \cong \mathbb{R}^{n+k}$. Then composing a given trivialisation $\tau_M \cong \mathbb{R}^n$ with $g \in O(n)$ corresponds to composing $\nu_M \cong \mathbb{R}^k$ with some $g' \in O(k)$ and a homotopy between
\[
\nu_j \oplus \tau_M \cong \mathbb{R}^{n+k} \xrightarrow{\ g \oplus g'\ } \mathbb{R}^{n+k}
\]
and the given trivialisation of $\nu_j \oplus \tau_M$, i.e.\ a nullhomotopy of $g \oplus g'$. We see that we need to use the inverse $O$-action. This distinction will be crucial below: with the right-action convention used in this paper, the tangential action is $\ominus J_n$, while the normal action is $J_n$.

\begin{remark}
\label{rmk:normal-action}
The subtlety about different $O$-actions is relevant since the map $\ominus \colon BO \to BO$ is not nullhomotopic. For example, consider the truncation $t \colon BO \to BO/B\mathrm{Spin} \cong B\mathbb{Z}/2 \times B\mathbb{Z}/2$ explicitly given by $(w_1, w_2)$. Then by the Whitney sum formula, composition of $\ominus$ with $t$ results in the map $(w_1, w_2 + w_1^2)$.
\end{remark}

\subsection{The abstract Pontryagin--Thom collapse map}
\label{sec:abstract-PT}

We use the abstract formulation of the Pontryagin--Thom collapse map in terms of Spanier--Whitehead duality, following \cite[Definition 10.3 and Proposition 10.5]{MR3252967}. Let $S$ denote the sphere spectrum.

\begin{definition}
\label{def:PT-collapse}
Let $X$ be a compact $d$-manifold, and fix an embedding $j \colon X \hookrightarrow \mathbb{R}^{d+k}$ with tubular neighbourhood $N_X$ and normal bundle $\nu_j$. The \emph{classical Pontryagin--Thom collapse map} is the based map
\[
c_j \colon S^{d+k} \longrightarrow \mathrm{Th}(\nu_j)
\]
obtained by collapsing the complement of $N_X$ in $\mathbb{R}^{d+k}$ (compactified to $S^{d+k}$) to the basepoint of the Thom space $\mathrm{Th}(\nu_j)$.

Passing to suspension spectra and using the canonical identification $\mathrm{Th}(\nu_j) \simeq \Sigma^{-d}(X^{-TX})$ (which holds because $\nu_j \oplus \tau_X \cong \mathbb{R}^{d+k}$), this gives the \emph{Pontryagin--Thom collapse map}
\[
c_X \colon S \longrightarrow X^{-TX},
\]
where $X^{-TX}$ denotes the Thom spectrum of the virtual bundle $-TX$ (the negative tangent bundle of $X$), i.e.\ the colimit
\[
X^{-TX} := \colim\!\left(X \xrightarrow{\ -TX\ } BO \times \mathbb{Z} \xrightarrow{\ J\ } \mathrm{Pic}(S) \to \mathrm{Sp}\right).
\]
The homotopy class of $c_X$ is independent of the choice of embedding $j$ and tubular neighbourhood $N_X$: any two choices are related by an isotopy, and the space of tubular neighbourhoods of a fixed embedding is contractible.
\end{definition}

\begin{proposition}
\label{prop:framed-PT}
Let $M$ be a closed $d$-manifold with a framing $\phi \colon TM \cong \mathbb{R}^d$. Then:
\begin{enumerate}
\item[$\natural_1$] The framing $\phi$ induces an equivalence $M^{-TM} \cong \Sigma^{-d}\Sigma^\infty_+ M$.
\item[$\natural_2$] 
The composition $S \xrightarrow{c_M} M^{-TM} \xrightarrow{\phi} \Sigma^{-d}\Sigma^\infty_+ M \xrightarrow{\text{crush}} \Sigma^{-d} S$ defines an element of $\pi_d S$. 
\item[$\natural_3$] This assignment is the classical Pontryagin--Thom isomorphism: it identifies the framed bordism group $\Omega_d^{\mathrm{fr}}$ with $\pi_d S$.
\end{enumerate}
\end{proposition}

\begin{proof}
$(\natural_1)$ follows from the fact that $-TM \cong -\mathbb{R}^d$ under the framing $\phi$, and the Thom spectrum of a trivial bundle of rank $-d$ is $\Sigma^{-d} \Sigma^\infty_+ M$. $\natural_2$ is the definition of the map: $c_M$ is the abstract PT collapse map of Definition~\ref{def:PT-collapse}, $\phi$ is the trivialisation induced by the framing, and the crush map is induced by $M \to \mathrm{pt}$. The fact that it is an isomorphism is the classical Pontryagin--Thom theorem \cite[Chapter IX]{Kosinski93}.
\end{proof}

\begin{lemma}
\label{lem:PT-equivariance}
For $g \in O(d)$, let $g \cdot \phi := \phi \circ g$ denote the framing obtained by acting on $\phi$ by the right action of $O(d)$. Then the map
\[
S \to \Sigma^{-d} S
\]
constructed in Proposition~\ref{prop:framed-PT} for the framing $g \cdot \phi$ is obtained from the map for $\phi$ by post-composing with the one-point compactification of $g^{-1} \colon \mathbb{R}^d \to \mathbb{R}^d$ (acting on the target $\Sigma^{-d} S = S^{-d}$). The resulting element of $\pi_d S$ is $\ominus J_d(g) := J_d(g)^{-1}$.
\end{lemma}

\begin{proof}
Changing the framing by $g$ changes the induced trivialisation of $-TM$ by $g^{-1}$. After the PT construction and crushing $M_+ \to S^0$, the induced map on the target $\Sigma^{-d} S = S^{-d}$ is the one-point compactification of $g^{-1} \colon \mathbb{R}^d \to \mathbb{R}^d$, whose stable class is $J_d(g)^{-1} = \ominus J_d(g)$.
\end{proof}

\begin{theorem}
\label{thm:PT-equivalence}
The Pontryagin--Thom construction induces an equivalence of spectra
\[
\mathrm{PT} \colon L(\mathrm{Bord}_n^{\mathrm{fr}}) \xrightarrow{\ \sim\ } S.
\]
\end{theorem}

\begin{proof}
The assignment $(M, \phi) \mapsto c_{M, \phi}$ of Proposition~\ref{prop:framed-PT} is compatible with disjoint union, so it defines a symmetric monoidal functor $\mathrm{PT} \colon \mathrm{Bord}_n^{\mathrm{fr}} \to \iota(S)$. Passing to localisations gives a map $\mathrm{PT} \colon L(\mathrm{Bord}_n^{\mathrm{fr}}) \to S$.

\emph{Step 1: Isomorphism on $\pi_0$.} The positively framed point maps to $1 \in \pi_0 S$.

\emph{Step 2: Isomorphism on $\pi_k$ for $k < n-1$.} By Corollary~\ref{cor:bordism-groups}, $\pi_k L(\mathrm{Bord}_n^{\mathrm{fr}}) \cong \Omega_k^{\mathrm{fr}}$. The classical Pontryagin--Thom theorem identifies $\Omega_k^{\mathrm{fr}} \cong \pi_k S$. For $k < n-1$, an $n$-framing on a $k$-manifold is equivalent to a stable framing, so the unstable and stable PT constructions agree on $\pi_k$.

\emph{Step 3: Isomorphism on $\pi_k$ for $k \geq n-1$.} Pass to the colimit over $n$. By Lemma~\ref{lem:bordism-colimit}, $\pi_k L(\mathrm{Bord}_\infty^{\mathrm{fr}}) \cong \Omega_k^{\mathrm{fr}} \cong \pi_k S$.

\emph{Step 4: Equivalence.} Since both sides are connective spectra and PT is an isomorphism on all homotopy groups, it is an equivalence.
\end{proof}
\subsection{Pontryagin--Thom: identifying \texorpdfstring{$L(\Bord_n^{\fr})$}{L(Bord)} with bordism groups}

The Pontryagin--Thom theorem gives a geometric interpretation of the homotopy groups of the sphere spectrum in terms of bordism groups of \emph{stably} framed manifolds.
On the other hand, $L(\Bord_n^{\fr})$ gives a geometric interpretation in terms of a certain spectrum built out of \emph{unstably} framed manifolds.
The goal of this section is to compare these descriptions and see how the Pontryagin--Thom theorem applies.

    Given a symmetric monoidal $(\infty,n)$-category $\mathcal{C}$, define the $(\infty,n-1)$-category $\Omega \mathcal{C} := \End_{\mathcal{C}} 1$.
    There is a map
     $\Omega^k \mathcal{C} \to \Omega^k L(\mathcal{C})$ given by $\Omega^k$ applied to the localisation map. %\luuk{I believe that I need to take $\Omega$ within connective spectra so strictly speaking taking $\pi_{\geq 0} \Omega^k X$}
     By the universal property of $L$ this gives a map $L(\Omega^k\mathcal{C}) \to \Omega^k L(\mathcal{C})$.
     For a general category, this map need not be close to an equivalence.
     However, if $\mathcal{C}$ has duals, then localisation identifies those as inverses, and consequently we do get:

\begin{lemma}
    If $k<n-1$, the map $L(\Omega^k\mathcal{C}) \to \Omega^k L(\mathcal{C})$ is an isomorphism on $\pi_0$.
    Explicitly, we obtain an isomorphism of groups $\pi_k L \mathcal{C} = \Omega^k \mathcal{C}/\sim$, where $x \sim y$ if and only if $x$ and $y$ are connected by a $k+1$-morphism.
\end{lemma}
\begin{proof}
    Since $k<n$, we can use the existence of adjoints of $k$-morphisms to ensure that $x,y \in \Omega^k \mathcal{C}$ as opposed to a zigzag of $k$-morphisms.
    There is no need to group-complete since $ev\colon x^\vee \otimes x \to 1$ will ensure $x^\vee$ is the inverse of $x$.
    Observe that $\sim$ is already an equivalence relation for $k < n-1$.
\end{proof}

For $\mathcal{C} = \Bord_n^{\fr}$, this means explicitly that $\pi_k L \Bord_n^{\fr}$ for $k<n$ are given by bordism groups.
Indeed, unpacking the definition of $\Omega^k \Bord_n^{\fr}$ we see that it has objects $k$-dimensional closed manifolds $Y$ with a trivialisation $TY \oplus \underline{\R}^{n-k} \cong \underline{\R}^n$, $1$-morphisms $k+1$-dimensional bordisms $M$ with a trivialisation $TM \oplus \underline{\R}^{n-k-1} \cong \underline{\R}^n$, etcetera.

\begin{definition}
    Let $B \to BO(n)$ be a tangential structure and $k<n$.
    The bordism group of $k$-dimensional manifolds with $B$-structure is the monoid of closed manifolds with $B$-structure under disjoint union modulo $k+1$-dimensional bordisms with $B$-structure.
\end{definition}

\begin{remark}
    $\sim$ need not be an equivalence relation for $k = n-1$.
    For example, consider the case $n = 2$ in $\Bord_2^{\fr}$.
    The manifold with a circle boundary obtained by puncturing a surface of genus $g\geq 2$ has an unstable framing, since it is homotopy equivalent to a $1$-dimensional CW-complex. We obtain a $1$-morphism in $\Bord_2^{\fr}$ from $S^1$ with some framing of $TS^1 \oplus \underline{\R}$ to $\emptyset$.
    There cannot be a $1$-morphism from $\emptyset$ to $S^1$ with the same framing, as their composition would result in an unstably framed closed manifold of genus $g' \geq g \geq 2$.
    
    However, for $k<n-1$, adjoints also exist for $k$-morphisms.
    Therefore, a morphism $x \to y$ exists if and only if a morphism $y \to x$ exists, and so $\sim$ is an equivalence relation.
    In particular, 
\end{remark}

\begin{corollary}\label{cor:bordism-groups}

If $k<n$, then $\pi_k L(\Bord_n^{\fr})$ is isomorphic to the bordism group of unstably framed manifolds.
If $k<n-1$, it is isomorphic to the bordism group of stably framed manifolds.
\end{corollary}
\begin{proof}
    The second point follows because in that case the $k+1$-morphisms come with a once stabilised trivialisation of the tangent bundle, which is equivalent to a stable trivialisation.
\end{proof}

Since we only understand $\pi_k L(\Bord_n^{\fr})$ in terms of bordism groups for small $k$, we take $n \to \infty$.
We comment briefly on the more mysterious geometric interpretations of $\pi_k L(\Bord_n^{\fr})$ for $k \geq n$.

\begin{example}
    For $n = 0$, the symmetric monoidal groupoid $\Bord_n^{\fr}$ of zero-dimensional manifolds and diffeomorphisms is equivalent to the groupoid 
    \[
    \bigsqcup_{k\geq 0} BS_k
    \]
    of finite sets equipped with a disjoint union.
    The fact that 
    \[
    L(\Bord_0^{\fr}) \cong S
    \]
    is a version of the Barratt-Priddy-Quillen theorem. 
\end{example}

\begin{example}
    The group $\pi_n L(\Bord_n^{\fr})$ is generated by zigzags of $n$-dimensional framed bordisms
    \[
    \emptyset \xrightarrow{M_1} Y_1 \xleftarrow{M_2} Y_2 \xrightarrow{M_3} \dots \xleftarrow{M_{k}} Y_{k} \xrightarrow{M_{k+1}} \emptyset
    \]
    between $(n-1)$-dimensional compact manifolds with a trivialisation $TY_j \oplus \underline{\R} \cong \underline{\R}^n$.
    The isomorphism of groups $\pi_n L(\Bord_n^{\fr}) \to \pi_n L(\Bord_{n+1}^{\fr}) \cong \Omega_n^{\fr}$ is explicitly given by taking the composition to a stably framed bordism $\emptyset \to \emptyset$ after potentially reversing bordisms by taking adjoints.
\end{example}

There is a symmetric monoidal functor $\Bord^{\fr}_n \to \Bord^{\fr}_{n+1}$ that maps a $k$-manifold $Y$ for $k\leq n$ with $TY \oplus \underline{\R}^{n-k} \cong \underline{\R}^n$ to the same manifold with the isomorphism $TY \oplus \underline{\R}^{n+1-k} \cong \underline{\R}^{n+1}$ obtained by stabilising with the identity, and a diffeomorphism of $n$-dimensional manifolds to its mapping cylinder as an $n+1$-dimensional manifold with the obvious unstable framing.
This functor is clearly $O(n)$-equivariant when restricting the $O(n+1)$-action on $\Bord^{\fr}_{n+1}$ to the first $n$ coordinates.
We thus get an $O$-action on
\[
\colim_n \Bord^{\fr}_n \in \Cat^{\E_\infty}_\infty
\]
where $\Cat^{\E_\infty}_\infty = \colim_n \Cat^{\E_\infty}_{(\infty,n)}$ is the $(\infty,1)$-category of $(\infty,\infty)$-categories defined as a colimit in $\Pr^L$.
We will be interested in the induced $O$-action on $L(\Bord^{\fr}_\infty) \cong \colim_n L(\Bord^{\fr}_n)$.

\begin{lemma}\label{lem:bordism-colimit}
    The degree $k$ homotopy group of the spectrum $L(\Bord^{\fr}_\infty)$ is given by $\Omega^{\fr}_k$, the bordism group of stably framed $k$-dimensional manifolds
\end{lemma}
\begin{proof}
    The functor $\pi_k \colon \mathcal{S}_* \to \Ab$ preserves filtered colimits because $S^k$ is compact, and so 
    \[
    \pi_k L(\colim_n \Bord^{\fr}_n) = \colim_n \pi_k L(\Bord^{\fr}_n).
    \]
    For $n$ large enough, $\pi_k L(\Bord^{\fr}_n) \cong \Omega^{\fr}_k$ and the structure maps in the colimit are the identity.
\end{proof}

In Section \ref{sec:abstract}, we concluded that there are equivalences $L(\Bord^{\fr}_n) \cong S$ for all $n$.
The map $\Bord^{\fr}_n \to \Bord^{\fr}_{n+1}$ induces the identity map on $S$ after applying $L$.
Indeed, since the cobordism hypothesis is given by evaluating on a point, the diagram
\[
\begin{tikzcd}
    \Fun(\Bord^{\fr}_n,X) \ar[r] \ar[d] & \Fun(\Bord^{\fr}_{n+1}, X) \ar[dl]
    \\
    X &
\end{tikzcd}
\]
commutes.
We therefore also get an equivalence $\colim_n L(\Bord_n^{\fr}) \cong S$.

Recall from Definition~\ref{def:PT-collapse} that the Pontryagin--Thom collapse map sends a closed normally framed $d$-manifold $Y \hookrightarrow \mathbb{R}^N$ to a based map $S^N \to S^{N-d}$, whose stable homotopy class lies in $\pi_d S$.

\begin{theorem}
    On homotopy groups, the equivalence $L(\Bord_\infty^{\fr}) \cong S$ is given by the Pontryagin--Thom construction.
\end{theorem}
\begin{proof}(Sketch)
    It suffices to show that the Pontryagin--Thom construction defines an equivalence of spectra $L(\Bord_\infty^{\fr}) \cong S$
    Indeed, we know that $S$ has only one nontrivial automorphism, and the Pontryagin--Thom construction sends $+$ to $1 \in \pi_0 S$.
     We therefore want to prove that the Pontryagin--Thom construction defines compatible $n$-dimensional TQFTs $\Bord_n^{\fr} \cong S$.
     
     We sketch a scanning argument for a model of $\Bord_n^{B}$ in which bordisms are embedded in $\R^\infty$ with an $n$-dimensional hypercube picked out for composition, and another infinitely many coordinates for the $\E_\infty$-structure.
     We will take $B$ to be the (stable) tangential structure of a (stable) normal framing; $\Bord_n^{\fr}$ clearly maps to this, and this is all we need for the PT construction.
     
     So suppose $X$ is an $n$-dimensional bordism embedded in $\R^N \times [0,1]$ from an $(n-1)$-dimensional closed manifold $Y_0$ embedded in $\R^N \times \{0\}$ to an $(n-1)$-dimensional closed manifold $Y_1$ embedded in $\R^N \times \{1\}$.
     The Thom-Pontryagin construction then gives a based map $S^N\times[0,1]\to S^{N-n}$ which restricts to the Thom-Pontryagin construction of $Y_1$ and $Y_2$ on the boundaries, which is the desired morphism in $[S^N,S^{N-n}]$.
     Observe that bordisms compose as homotopies do.
     A similar construction holds for an arbitrary $k$-morphism contained in $[0,1]^k \times \R^N$ giving a map $S^N\times[0,1]^k\to S^{N-n-k}$.
     The space of tubular neighbourhoods is contractible, so the construction is well-defined up to coherent homotopy.
\end{proof}

\subsection{Identifying the \texorpdfstring{$O(n)$}{O(n)}-action}

\begin{theorem}[Main theorem]
\label{thm:main}
With respect to the right action of $O(n)$ on framings (Convention~\ref{conventions-item:frame}), the $O(n)$-action on $L(\mathrm{Bord}_n^{\mathrm{fr}}) \simeq S$ induced by change of tangential framing is given by
\[
\ominus J_n \colon O(n) \longrightarrow \mathrm{GL}_1(S), \qquad g \mapsto J_n(g)^{-1}.
\]
Passing to the colimit as $n \to \infty$, the $O$-action on $\mathrm{colim}_n L(\mathrm{Bord}_n^{\mathrm{fr}}) \simeq S$ is given by $\ominus J \colon O \to \mathrm{GL}_1(S)$.
\end{theorem}

\begin{proof}
By Theorem~\ref{thm:PT-equivalence}, the Pontryagin--Thom construction gives an equivalence of spectra $\mathrm{PT} \colon L(\mathrm{Bord}_n^{\mathrm{fr}}) \to S$, under which the positively framed point corresponds to $1 \in \pi_0 S$. By Lemma~\ref{lem:PT-equivariance}, this equivalence is $O(n)$-equivariant, with the action of $g$ on $S$ given by $\ominus J_n(g) := J_n(g)^{-1}$.

More explicitly, on the standard $n$-ball $B^n \subset \mathbb{R}^n$ with its standard framing, the PT map is the identity $S^n \to S^n$ (as an element of $\pi_n S$). Changing the framing by $g$ (right action) rotates the target sphere by $g^{-1}$, giving $\ominus J_n(g) \in \mathrm{GL}_1(S)$. Since the $n$-balls generate $S$ as a spectrum, the action on all of $S$ is $\ominus J_n(g)$.

Passing to the colimit as $n \to \infty$ gives the stable $J$-homomorphism composed with $\ominus$: $\ominus J \colon O \to \mathrm{GL}_1(S)$.
\end{proof}
\begin{remark}
\label{rmk:rotation-picture}
The $O(n)$-action can be seen directly in the classical Pontryagin--Thom construction. Embed $M^d \hookrightarrow \mathbb{R}^{d+k}$ with a normal framing. Changing the tangential framing of $M$ by $g \in O(d)$ is equivalent to rotating the embedding coordinates by $g \oplus \mathrm{id}_k \in O(d+k)$, which acts on the source sphere $S^{d+k}$ in the collapse map $S^{d+k} \to S^k$. After stabilisation, this rotation is the unstable $J$-homomorphism $J_d(g)$. With the right-action convention for framings, the normal framing changes by $g^{-1}$, so the collapse map is post-composed with $J_d(g)^{-1}$, giving $\ominus J_d(g)$ (Theorem~\ref{thm:main}).
\end{remark}
\begin{remark}
\label{rmk:Lurie-observation}
The idea that the $O(n)$-action on the space of fully dualisable objects in a symmetric monoidal $(\infty,n)$-category with duals is related to the $J$-homomorphism that appears in Lurie's original treatment of the cobordism hypothesis \cite[Example 2.4.15]{MR2555928}. Lurie sketches that for a Picard $\infty$-groupoid, the $O(n)$-action is ``a version of the classical $J$-homomorphism''. In this paper we make this identification explicit: using the right-action convention for the $O(n)$-action on framings, we identify the action as $\ominus J_n$.
\end{remark}

\begin{remark}
\label{rmk:GP-comparison}
An equivalent result has been obtained by Grady--Pavlov \cite[Theorem 5.0.6]{grady2021geometric} using the framework of geometric cobordism categories. The proof above works directly within Lurie's formulation of the cobordism hypothesis, using the abstract Pontryagin--Thom collapse map. Our convention (right action of $O(n)$ on framings) gives the action $\ominus J_n$; in the left-action convention, the same proof gives $J_n$.
\end{remark}
\subsection{The classification of invertible field theories}

An invertible TQFT is a TQFT in which the target has invertible $k$-morphisms for all $k \geq 0$. More precisely, an $n$-dimensional $B$-structured invertible TQFT with target a connective spectrum $X$ is a symmetric monoidal functor $\Bord^B_n \to \iota X$, or equivalently, a map of spectra $L\Bord^B_n \to X$.

The following result is at the cornerstone of modern computations of anomalies and invertible phases of matter \cite{freed2021reflection}, and was shown using different techniques by \cite{MR4821355}.

\begin{proposition}
\label{prop:MTB}
If $\xi\colon B \to BO(n)$ is a tangential structure, there is an equivalence of spectra $L\Bord^B_n \cong MTB$ natural in $B$.
\end{proposition}
\begin{proof}
We use the equivalence $\mathrm{Bord}_n^B \simeq \colim_B \mathrm{Bord}_n^{\mathrm{fr}}$. Since $L$ preserves colimits, we get
\[
L(\mathrm{Bord}_n^B) \simeq \colim_B L(\mathrm{Bord}_n^{\mathrm{fr}}) \simeq \colim_B S,
\]
where the $B$-action on $S$ is via the composition
\[
B \xrightarrow{\ \xi\ } BO(n) \xrightarrow{\ (\iota, n)\ } BO \times \mathbb{Z} \xrightarrow{\ \ominus\ } BO \times \mathbb{Z} \xrightarrow{\ J\ } \mathrm{Pic}(S) \longrightarrow \mathrm{Sp},
\]
using Theorem~\ref{thm:main} (with the right-action convention, the $O(n)$-action is $\ominus J_n$, and passing to $B$ gives the composition above). By Definition~\ref{def:MTB}, this colimit is exactly the Madsen--Tillmann spectrum $MTB$ (the Thom spectrum of the virtual normal bundle $-\xi$). Therefore
\[
L(\mathrm{Bord}_n^B) \simeq MTB.
\]
Naturality in $B$ follows from the naturality in each step.
\end{proof}

\begin{remark}
\label{rmk:MTB-convention}
For the classification of invertible TFTs, the identification $L(\mathrm{Bord}_n^B) \simeq MTB$ means that a $B$-structured invertible TFT is a map $MTB \to X$ of spectra, i.e, an element of $\pi_0 \mathrm{Map}_{\mathrm{Sp}}(MTB, X)$. Equivalently, it is a homotopy fixed point $X^{hB}$, where $B$ acts on $X$ via $\ominus J_n$ (Theorem~\ref{thm:iTFT-classification}).
\end{remark}
\begin{remark}
\label{rmk:J-vs-ominusJ}
With the right-action convention for the $O(n)$-action on framings (Convention~\ref{conventions-item:frame}), the $O(n)$-action on $L(\mathrm{Bord}_n^{\mathrm{fr}})$ is $\ominus J_n$, not $J_n$ (Theorem~\ref{thm:main}). This is because the change of framing acts by $g^{-1}$ on the induced trivialisation of the tangent bundle. Equivalently, one can use the \emph{left}-action convention for the $O(n)$-action on framings, in which case the action is $J_n$; the two conventions differ by the inversion automorphism $\ominus$ (Convention~\ref{conventions-item:O-op}). In this paper we use the right-action convention, so the action is $\ominus J_n$.
\end{remark}
Now, it follows that 

\begin{theorem}
\label{thm:iTFT-classification}
The spectrum of $B$-structured invertible topological field theories valued in a connective spectrum $X$ is equivalent to
\[
\mathrm{Map}_{\mathrm{Sp}}(L(\mathrm{Bord}_n^B), X) \simeq \mathrm{Map}_{\mathrm{Sp}}(MTB, X),
\]
where the $B$-action on $X$ is via the composition
\[
B \xrightarrow{\ \xi\ } BO(n) \xrightarrow{\ B(\ominus J_n)\ } B\mathrm{GL}_1(S) \xrightarrow{\ \otimes X\ } \mathrm{Sp}.
\]
Equivalently, the spectrum is the homotopy fixed point spectrum $X^{hB}$.
\end{theorem}

\begin{proof}
By Proposition~\ref{prop:MTB}, $L(\mathrm{Bord}_n^B) \simeq MTB$. The spectrum of $B$-structured invertible TFTs valued in $X$ is therefore
\[
\mathrm{Map}_{\mathrm{Sp}}(L(\mathrm{Bord}_n^B), X) \simeq \mathrm{Map}_{\mathrm{Sp}}(MTB, X).
\]
For $B = \mathrm{pt}$, we have $L(\mathrm{Bord}_n^{\mathrm{fr}}) \simeq S$, so the spectrum is $\mathrm{Map}_{\mathrm{Sp}}(S, X) \simeq X$. The $O(n)$-action is induced by precomposition with the $\ominus J_n$-action on $S$. By adjunction, this is the action on $X$ via the $S$-module structure. Taking homotopy fixed points yields $X^{hB}$.
\end{proof}
\begin{remark}[Conjectural]
\label{rem:Pic-n-S}
Given a connective commutative ring spectrum $R$, define $\mathrm{Mod}^{(n)}(R)$ inductively as the higher idempotent completion \cite{gaiotto2019condensations} of the delooping of the symmetric monoidal $(\infty,n-1)$-category $\mathrm{Mod}^{(n-1)}(R)$, where $\mathrm{Mod}^{(0)}(R) := R$. Let $\mathrm{Pic}^{(n)}(R)$ be the maximal Picard sub-$\infty$-groupoid of $\mathrm{Mod}^{(n)}(R)$. Then $\mathrm{Pic}^{(2)}(R)$ is the Brauer--Picard spectrum.

For $R = S$, $\mathrm{Pic}^{(2)}(S)$ is known to be $B\mathrm{Pic}(S)$ \cite{MR3190610}. As far as we know, $\mathrm{Pic}^{(n)}(S)$ is not known for $n > 2$. The $O(n)$-action on this space will, by the cobordism hypothesis, be induced from the $O(n)$-action on $\mathrm{Bord}_n^{\mathrm{fr}}$, whose localisation is computed by Theorem~\ref{thm:main} to be $\ominus J_n$.
\end{remark}
\section{Computations}\label{sec:computations}

In this section, we illustrate Theorem~\ref{thm:main} by determining the $O(n)$-action on the homotopy groups of the sphere spectrum $S$ in low degrees. By Theorem~\ref{thm:main}, the $O(n)$-action on $L(\Bord_n^{\fr}) \simeq S$ is given by the composition
\[
O(n) \xrightarrow{\ \ominus J_n\ } \GL_1(S) \longrightarrow \Aut(S),
\]
where the second map sends a unit to the corresponding multiplication map on $S$.
 We determine the effect of this action on $\pi_k S$ for small $k$ and all $n$.

\subsection{The J-homomorphism concretely}

The sphere spectrum $S$ plays the role of the integers in abelian groups; it is the initial $\E_\infty$-ring spectrum, and every spectrum $X$ has a canonical $S$-module structure.
If $X \in \Sp$, the map $S \times X \to X$ induces action maps $\pi_{k} S \times \pi_{l} X \to \pi_{k+l} S$.

Recall the homotopy groups of the sphere spectrum in low degrees~\cite{ravenel2003complex}
\[
\pi_0 S \cong \mathbb{Z}, \qquad
\pi_1 S \cong \mathbb{Z}/2 \;\; (\text{generated by } \eta), \qquad
\pi_2 S \cong \mathbb{Z}/2 \;\; (\text{generated by } \eta^2), \qquad
\pi_3 S \cong \mathbb{Z}/24 \;\; (\text{generated by } \nu).
\]

\begin{lemma}
For $k = 0$, the action map is given by $-1 \colon \pi_l X \to \pi_l X$.
\end{lemma}

For $k>0$, the action is much more mysterious. 
For example, for for $k = 1$ the generator of $\pi_1 S\cong \Z/2$ gives a map
    $\pi_l X \to \pi_{l+1} X$ which by bilinearity of the action map factors through a map $\eta_{X,l} \colon \pi_l X \to \pi_{l+1} X/2$.

\begin{lemma}
\label{lem:eta}
The stable cohomology group $\pi_0 Map_{Sp}(\pi_l X, \Sigma \pi_{l+1} X)$ is isomorphic to $\Hom_{Ab}(\pi_l X, \pi_{l+1} X/2)$.
The element $\eta_{X,l}$ of this group is the $k$-invariant of $X$ at the spot $l,l+1$.
\end{lemma}

We are not aware of a general way to describe the action $\pi_l X \to \pi_{l+3} X/24$ of $\nu$.

Definition \ref{def:units} gives the pullback
\[
\begin{tikzcd}
\GL_1(S) \ar[r] \ar[d] & \Omega^\infty S \ar[d, "\pi_0"] \\
\{\pm 1\} \ar[r] & \mathbb{Z}
\end{tikzcd}
\]
in the $\infty$-category $\mathcal{S}$. From the associated fibre sequence, one reads off its homotopy groups:
\[
\pi_0 \GL_1(S) \cong \mathbb{Z}/2 \cong \{\pm 1\}, \qquad
\pi_k \GL_1(S) \cong \pi_k S \quad \text{for all } k \ge 1.
\]
The image of the $J$-homomorphism is understood better than the homotopy groups of $\GL_1(S)$. 
In low degrees, we have \cite{adams1963groups} 
\begin{align*}
\Z \cong \pi_0 (BO \times \Z) \xrightarrow{\cong}& \pi_0 \Pic(S) \cong \Z
\\
    \Z/2 \cong \pi_0 O \xrightarrow{\cong}& \pi_0 GL_1(S) \cong \Z/2
    \\
    \Z/2 \cong \pi_1 O \xrightarrow{\cong}& \pi_1 GL_1(S) \cong \Z/2 \ni \eta
    \\
    0 \cong \pi_2 O \xrightarrow{0}& \pi_2 GL_1(S) \cong \Z/2 \ni \eta^2
    \\
    \Z \cong \pi_3 O \twoheadrightarrow & \pi_3 GL_1(S) \cong \Z/24 \ni \nu
    \\
    0 \cong \pi_4 O \xrightarrow{0}& \pi_4 GL_1(S) \cong 0
    \\
    0 \cong \pi_5 O \xrightarrow{0}& \pi_5 GL_1(S) \cong 0
    \\
    0 \cong \pi_6 O \xrightarrow{0}& \pi_6 GL_1(S) \cong \Z/2 \ni \nu^2
    \\
    \Z \cong \pi_7 O \twoheadrightarrow& \pi_7 GL_1(S) \cong \Z/240
\end{align*}

We obtain:

\begin{proposition}
\label{prop:action-summary}
\begin{enumerate}
    \item[$\natural_1$] The generator of $\pi_0 O$ acts on $X$ by $-1$.
    \item[$\natural_2$] The generator of $\pi_1 O$ acts on $X$ by $\eta$ as described in Lemma \ref{lem:eta}.
    \item[$\natural_3$] The generator of $\pi_3 O$ acts on $X$ by $\nu$.
\end{enumerate}
\end{proposition}
The full $O$-action on $X$ contains even more interesting information about how the homotopy groups of $BO$ fit together.

\begin{remark}\label{rem:spectral-sequence}
The computation in Proposition~\ref{prop:action-summary} determines the $O(n)$-action on each $\pi_k S$ as a self-map. The full homotopy fixed point spectrum $S^{hB}$ for a tangential structure $B \to BO(n)$ is computed by the homotopy fixed point spectral sequence
\[
E_2^{p,q} = H^p(B; \pi_{-q} S) \implies \pi_{-p-q}(S^{hB}),
\]
where the coefficients $\pi_* S$ carry the $B$-action described above. This spectral sequence computes the homotopy groups of the spectrum of $B$-structured invertible TQFTs valued in $S$. The spectral sequence also involves higher coherences of the action: the fact that $B$ acts via the map $B \to BO(n) \xrightarrow{\ominus J_n} \GL_1(S)$ means that elements of $\pi_1(B)$ can act by multiplication by $\eta \in \pi_1 S$, and elements of $\pi_3(B)$ by $\nu \in \pi_3 S$. These enter the spectral sequence through differentials and Massey products, but do not affect the action on individual homotopy groups as self-maps.

For $n=2$, a rotation with odd winding number has image $\eta \in \pi_1 \GL_1(S)$. For $n \ge 3$, the nontrivial element of $\pi_1(SO(n))$ maps to $\eta$ \cite{adams1963groups}. In both cases, this $\pi_1$-information contributes to the $d_2$-differential in the spectral sequence, encoding the fact that a $2\pi$ rotation acts nontrivially on the \emph{space} of invertible TFTs, even though it acts trivially on the homotopy groups of $S$ as self-maps. This is the mathematical mechanism behind the spin-statistics theorem in the context of framed field theories \cite{freed2021reflection}.

There is a direct relationship between this spectral sequence and the Madsen--Tillmann spectrum $MTB$. By the classification of invertible TFTs (Theorem~\ref{thm:iTFT-classification}), we have $L(\mathrm{Bord}_n^B) \simeq MTB$, and therefore
\[
S^{hB} \simeq \mathrm{Map}_{\mathrm{Sp}}(MTB, S).
\]
Consequently, the homotopy fixed point spectral sequence converges to the homotopy groups of the function spectrum $\mathrm{Map}_{\mathrm{Sp}}(MTB, S)$. For specific choices of $B$, the differentials in this spectral sequence encode the bordism groups of $B$-structured manifolds: the $d_2$-differential, governed by the $\eta$-action, obstructs the existence of certain $B$-structured invertible TFTs in dimension $n$; higher differentials involving $\nu$ obstruct the extension of such theories to dimension $n+1$.
\end{remark}

\subsection{What these computations determine}

In the previous subsections we computed the action for target $S$; we now explain how this generalises to an arbitrary connective spectrum $X$.
By Corollary~\ref{cor:CH-tangential} applied to $\mathcal{C} = \iota(X)$ (so that $\mathcal{C}^{\fd} \simeq \Omega^\infty X$ and the homotopy fixed points are taken at the spectrum level), for a tangential structure $\xi\colon B \to BO(n)$, the space of $B$-structured invertible TFTs valued in a spectrum $X$ is $X^{hB}$, where $B$ acts via $\xi$ followed by $\ominus J_n$. Proposition~\ref{prop:action-summary} determines this action on $\pi_k {S}$ for $k \le 3$. For $X = {S}$, this provides the coefficient module for the homotopy fixed point spectral sequence computing ${S}^{hB}$ in low degrees.

\section{Physical Interpretations}\label{sec:physics}
In this section, we briefly indicate how our main theorem (Theorem~\ref{thm:main}) connects to the physics of invertible topological phases and anomalies. The discussion is intentionally informal and aims to translate our mathematical results into the language of field theory; we make no claim to physical rigour. Conjectural statements are explicitly marked.

\subsection{Invertible field theories and anomalies}

An invertible $n$-dimensional topological field theory is one where the state spaces are invertible under the tensor product. In physics, invertible TFTs describe the low-energy limits of symmetry-protected topological (SPT) phases and the anomalies of quantum field theories; see \cite{KapustinThorngren, Debray, BeaudryCampbell, Freed} for the physics. Mathematically, they are classified by the cobordism hypothesis: for a tangential structure $B \to BO(n)$, the space of $B$-structured invertible TFTs valued in a connective spectrum $X$ is $X^{hB}$, where $B$ acts on $X$ via $\xi$ followed by $\ominus J_n$ (with our convention). In the lowest dimensions, the classification becomes concrete, and we unpack it below.

\subsection{Spacetime dimension \texorpdfstring{$1$}{1}: time-reversal and fermion parity}

For $n=1$, $O(1)\cong \mathbb{Z}/2$, generated by a reflection. By Proposition~\ref{prop:action-summary}, the action on $\pi_k S$ is multiplication by $-1$:

\begin{itemize}
    \item[$\natural_1$] $\pi_0 S \cong \mathbb{Z}$: $x \mapsto-x$, fixed subgroup $\{0\}$.
    \item[$\natural_2$] $\pi_1 S \cong \mathbb{Z}/2$: trivial, fixed subgroup $\mathbb{Z}/2$.
    \item[$\natural_3$] $\pi_2 S \cong \mathbb{Z}/2$: trivial, fixed subgroup $\mathbb{Z}/2$.
    \item[$\natural_4$] $\pi_3 S \cong \mathbb{Z}/24$: $x \mapsto -x$, fixed subgroup $\mathbb{Z}/2$.
\end{itemize}

Physically, $O(1)$ is time-reversal symmetry. The invariant on $\pi_0$ classifies the gravitational response: theories with nonzero $\pi_0$-component are not time-reversal symmetric, reflecting the fact that time-reversal complex conjugates the partition function on a closed $1$-manifold. The $\mathbb{Z}/2$ invariant on $\pi_1$ (the $\eta$-invariant) is time-reversal invariant, consistent with the fact that the Kitaev chain and other $1$d fermionic SPT phases are protected by time-reversal. The $\mathbb{Z}/24$ invariant on $\pi_3$ has a $\mathbb{Z}/2$ subgroup fixed by time-reversal, corresponding to the bosonic SPT phases that survive time-reversal symmetry.

\begin{remark}
The fact that the $\pi_1 \cong \mathbb{Z}/2$ invariant is fixed by time-reversal matches the physics of the Kitaev chain: the boundary of the nontrivial phase is a $0+1$d system whose effective description is captured by this invariant. Our computation takes place in the fully framed setting, which captures the gravitational (framing) anomaly sector. A complete treatment of fermionic SPT phases would replace the framed bordism category with a spin or $\mathrm{pin}^\pm$ bordism category; our main theorem applies to any tangential structure $B \to BO(n)$ via Corollary~\ref{cor:CH-tangential}, so such refinements are accessible within the framework of this paper.
\end{remark}
\subsection{Spacetime dimension \texorpdfstring{$2$}{2}: framing anomaly and the spin statistics connection}

For $n=2$, $O(2)$ has $\pi_0 \cong \mathbb{Z}/2$. Rotations act trivially on all $\pi_k S$; reflections act by multiplication by $-1$, identical to the $n=1$ case. The action on low homotopy groups is therefore the same as above.

The physically relevant invariants for $2$-dimensional framed theories are:
\begin{itemize}
    \item[$\natural_1$] $\pi_0 S \cong \mathbb{Z}$: a reflection-invariant theory must have trivial gravitational response (the partition function on a closed $2$-manifold is invariant under orientation reversal).
    \item[$\natural_2$] $\pi_1 S \cong \mathbb{Z}/2$: the $\eta$-invariant classifies $2$d fermionic invertible phases; it is reflection-invariant.
    \item[$\natural_3$] $\pi_3 S \cong \mathbb{Z}/24$: this is the \emph{framing anomaly} (or chiral central charge mod $24$). The reflection-invariant phases are those with framing anomaly $0$ or $12$ mod $24$. 
\end{itemize}

\begin{remark}
    Although rotations act trivially on each $\pi_k S$ as a self-map, the higher homotopy coherence of the $O(2)$-action (Remark~\ref{rem:spectral-sequence}) means that a $2\pi$ rotation acts by $\eta \in \pi_1 S$ in the homotopy fixed point spectral sequence. This is the mathematical reflection of the spin-statistics theorem: a $2 \pi$ rotation can act non-trivially on the \emph{space} of states, even though it acts trivially on the partition function (See Remark~\ref{rem:framing-anomaly} below). For $n=2$, this $\eta$-action contributes to the $d_2$ differential and can obstruct the existence of rotation-equivariant theories.
\end{remark}

\begin{remark}
\label{rem:framing-anomaly}
The framing anomaly of a $2$d chiral theory is characterised by the failure of the theory to be defined on an arbitrary framed manifold without a choice of $3$-dimensional bulk. As Witten explains \cite[pp.~191--192]{witten1992holomorphic}, the Hilbert space of a chiral CFT depends on a choice of rigging (framing) at the boundary: under a $2\pi$ rotation of the rigging, a state of conformal weight $h$ transforms as $|\psi\rangle \mapsto e^{2\pi i h}|\psi\rangle$, while the partition function on a closed surface is independent of the rigging.

In the framework of this paper, the anomaly of a $2$d framed TFT valued in $S$ is classified by $\pi_3 S \cong \mathbb{Z}/24$: the theory assigns an element of $\pi_3 S$ to a framed $3$-manifold via the TFT's extension to the point, and this element is precisely the framing anomaly (or chiral central charge modulo $24$). A $2$d theory can be realised as the boundary of a $3$d invertible theory if and only if its framing anomaly vanishes, recovering the classic anomaly inflow result \cite{witten1989quantum}. The $\mathbb{Z}/2$ fixed subgroup of the $O(2)$-action on $\pi_3 S$ (Proposition~\ref{prop:action-summary}) means that reflection-symmetric $2$d theories can have framing anomaly $0$ or $12$ modulo $24$, while theories with other values of the framing anomaly break orientation-reversal symmetry.
\end{remark}

\begin{remark}\label{rem:Z/24}
    Witten, in \cite[Equations 2.24, 2.25, 4.41]{witten1989quantum} shows that under a change of framing by $s$ units, the partition function transforms by the phase $\text{exp}(2\pi is\cdot c/24)$, so the framing anomaly is classified by $\mathbb{Z}/24$.
\end{remark}

\subsection{Spacetime dimension \texorpdfstring{$3$}{3} and beyond}

For $n \ge 3$, the $O(n)$-action on $\pi_k S$ is identical to the $n=2$ case: rotations act trivially, reflections act by multiplication by $-1$. The physical interpretation of the low-degree invariants extends straightforwardly from the $2$-dimension discussion.

The higher homotopy coherence becomes richer in dimension $3$: the $J$-homomorphism on $\pi_3(SO(n))$ (for $n \ge 3$) involves the element $\nu \in \pi_3 S \cong \mathbb{Z}/24$, which enters the homotopy fixed point spectral sequence through higher differentials \cite{adams1963groups}. This suggests that $3$-dimensional invertible phases with crystalline symmetries may exhibit more subtle mechanisms for anomaly cancellation.

\subsection{Relation to the Adams Conjecture and Chromatic Phenomena}

The $J$-homomorphism is intimately related to the \emph{Adams Conjecture} and the image of $J$ in the stable homotopy groups of spheres \cite{adams1963groups}. Our main theorem shows that the $O(n)$-action on the space of invertible TFTs valued in $S$ factors through the image of $\ominus J$. This suggests that the classification of SPT phases is controlled by the \emph{chromatic filtration} of the sphere spectrum: the image of $\ominus J$ captures the height-$1$ chromatic layer, while higher chromatic phenomena (detected by Morava $E$-theories) could classify more exotic invertible phases.

\bibliographystyle{alpha}
\bibliography{biblio.bib}{}

\end{document}